\documentclass[letterpaper]{amsart}
\usepackage{amssymb}
\usepackage{graphics}
\usepackage{pst-plot,pst-node}

\newtheorem{theorem}{Theorem}[section]
\newtheorem{lemma}[theorem]{Lemma}

\theoremstyle{definition}

\newtheorem{remark}{Remark}

\newtheorem{question}{Question}

\newcommand{\euO}{\mathfrak O}
\newcommand{\euP}{\mathfrak P}

\newcommand{\bQ}{\mathbb Q}
\newcommand{\bZ}{\mathbb Z}
\newcommand{\bF}{\mathbb F}

\begin{document}

\title[Quaternion Extensions]{On wild ramification in
quaternion extensions}

\author[G. Griffith {\sc Elder}]{{\sc G. Griffith} Elder}
\author[Jeffrey J. {\sc Hooper}]{{\sc Jeffrey J.} Hooper}

\thanks{Elder was partially supported by NSF grant DMS-0201080, Hooper
was partially supported by NSERC}


\address{G. Griffith {\sc Elder}\\
Department of Mathematics \\
Virginia Tech \\
Blacksburg, VA 24061-0123 U.S.A.}
\email{elder@vt.edu}


\address{Jeffrey J. {\sc Hooper}\\
Department of Mathematics and Statistics \\
Acadia University\\
Wolfville, NS \\
B4P 2R6\\
Canada}
\email{jeff.hooper@acadiau.ca}

\date{September 9, 2006}

\begin{abstract}
This paper provides a complete catalog of the break numbers that
occur in the ramification filtration of fully and thus wildly
ramified quaternion extensions of level one, dyadic number fields
(along with some partial results for level $>1$).  This catalog
depends upon the {\em refined ramification filtration}, which as
defined in \cite{elder:newbreaks} is associated with the biquadratic
subfield. Moreover we find that quaternion counter-examples to the
conclusion of the Hasse-Arf Theorem are extremely rare and can occur
only when the refined ramification filtration is, in two different
ways, extreme.
\end{abstract}

\maketitle

\section{Introduction}
Quaternion extensions are often the smallest extensions to exhibit
special properties and have played an important role in Galois module
structure \cite{frohlich}.  In the setting of the Hasse-Arf Theorem,
they are used to illustrate the fact that upper ramification numbers
in a non-abelian extension need not be integers \cite[IV\S4 Exercise
2]{serre:local}.  To better understand the counter-examples to the
conclusion of the Hasse-Arf Theorem and as a first step towards an
explicit description of wildly ramified Galois module structure
({\em e.g.} \cite{elder:onebreak,elder:annals,elder:biquad,elder:c8}),
we catalog the ramification break numbers of totally
ramified quaternion extensions of dyadic number fields.

\subsection{Notation}
Let $\bQ_2$ be the field of dyadic numbers, and let $K/\mathbb{Q}_2$
be a finite extension with $T$ its maximal unramified subfield. Then
$e_K=[K:T]$ is its degree of absolute ramification and
$f_K=[T:\mathbb{Q}_2]$ is its degree of inertia.  We will continue to
use subscripts to denote field of reference. So $\pi_K$ is a prime
element in $K$, $\euO_K$ the ring of integers, $\euP_K=\pi_K\euO_K$
its maximal ideal, and $v_K(\cdot)$ the valuation normalized so that
$v_K(\pi_K^n)=n$ for $n\in\bZ$. By abuse of notation, we identify the
residue fields $\euO_K/\euP_K=\euO_T/\euP_T=\bF_q$ with the finite
field of $q=2^{f_K}$ elements.

Let $N/K$ be a fully ramified quaternion extension with
$$G=\mbox{Gal}(N/K)=\left\langle\sigma,\gamma \mid \sigma^2=\gamma^2,
\gamma^{-1}\sigma\gamma=\sigma^{-1}\right\rangle.$$ It is a quick
exercise to check that these relations, $\sigma^2=\gamma^2$ and
$\gamma^{-1}\sigma\gamma=\sigma^{-1}$, yield $\sigma^4=1$.  Recall the
ramification filtration $G_i=\{s\in G:v_N((s-1)\pi_N)\geq i+1\}$ and
that {\em break} {\em numbers} (or jump numbers) are those integers
$b$ such that $G_b\supsetneq G_{b+1}$ \cite[ChIV]{serre:local}. .

Since $\mbox{Gal}(N/K)$ has a unique subgroup of order $2$, namely
$\langle\sigma^2\rangle$, and since the quotient of consecutive
ramification groups (in a fully ramified $p$-extension) is necessarily
elementary abelian \cite[IV \S2 Prop 7 Cor]{serre:local}, the
ramification filtration for $N/K$ decomposes naturally into two
filtrations: one for $M/K$ where $M=N^{\langle\sigma^2\rangle}$, and
one for $N/M$. Indeed the ramification break for $N/M$ is the largest
ramification break for $N/K$ \cite[IV \S1 Prop 2]{serre:local}. The
other break(s) for $N/K$ are those of $M/K$ \cite[IV \S1 Prop 3
Cor]{serre:local}.  This suggests

\begin{question}
How does ramification above ({\em i.e.} in $N/M$) depend upon
ramification below ({\em i.e.} in $M/K$)?
\end{question}

\subsection{On Ramification in Biquadratic Extensions}
There are either one or two break numbers in the ramification
filtration for the quotient group $\overline{G}=\mbox{Gal}(M/K)$.  In
the one break case, the break satisfies $1\leq b<2e$ with $b$ odd. In
the two break case, the breaks $b_1<b_2$ satisfy $1\leq b_1<2e$ with
$b_1$ odd, and $b_1<b_2\leq 4e-b_1$ with $b_2\equiv b_1\bmod 4$ when
$b_2< 4e-b_1$. This follows by considering upper ramification numbers and
the Herbrand Function \cite[IV\S3]{serre:local}.

Let $b_3$ denote the break for $\mbox{Gal}(N/M)$. Then from \S1.1 we
see that the ramification breaks for $G$ are either $b<b_3$ or
$b_1<b_2<b_3$.  To give a complete description of $b_3$ in the one
break case, $b<b_3$, we will need information provided by the refined
ramification filtration \cite{elder:newbreaks}. This is discussed in
detail as part of \S3, so for now we simply summarize the main
results: (1) There is a refined second break number $r\in \bZ$, which
satisfies $b<r<b_3$. (2) Associated with this second refined break
number is a $q-1$ root of unity $\omega$ (actually an equivalence
class, but for the moment it does no harm to confuse the equivalence class
with its representative).

As a result, to any fully ramified quaternion extensions of $N/K$ we
can assign a ramification triple: either $(b,r,b_3)$ in the one break
case or $(b_1, b_2, b_3)$ in the two break case and we are interested
in cataloging these triples. Our catalog has three cases depending
upon ramification in $M/K$.
\begin{itemize}
\item If $M/K$ has one ramification break $b$, then there is a second
refined break $r$ along with an associated root of unity $\omega$.
\begin{itemize}
\item If $\omega^3= 1$, we place $N$ in $\mathcal{Q}_{1^*}^K$.
\item If $\omega^3\neq 1$, we place $N$ in $\mathcal{Q}_1^K$.
\end{itemize}
\item If $M/K$ has two ramification breaks $b_1<b_2$, then
we place $N$ in $\mathcal{Q}_2^K$.
\end{itemize}
And so $\mathcal{Q}_{1^*}^K\cup\mathcal{Q}_1^K\cup \mathcal{Q}_2^K$ is the
set of all fully ramified quaternion extensions of $K$.

\subsection{Catalog of triples: Subsets of $\bZ^3$}\label{ssec:catalog}
In this section, based upon a choice of
positive integer $e$, we define three sets of triples
$\mathcal{R}_i^e\subset \bZ^3$ with $i\in \{1, 1^*, 2\}$, whose
elements will be denoted by $(s_1,s_2,s_3)$. In our descriptions of
these sets, the values of a coordinate $s_j$ will, in each case,
depend upon the values of preceding coordinates ($s_h$ for $h<j$).
So we begin by describing the first coordinates. In all cases
$$s_1\in S_1=\{n\in \bZ: 0< n<2e, n\equiv 1\bmod 2\}.$$ To
describe $s_2$, we must consider two basic cases: $i\in\{1^*, 1\}$ and $i=2$.
Let
$$m_i(s_1)=\begin{cases}
\min\{2s_1,4e-s_1\}&\mbox{for }i=1^*, 1,\\
4e-s_1&\mbox{for }i=2.\end{cases}$$
Then
$$s_2\in S_2^i(s_1)=\{n\in \bZ: s_1< n\leq m_i(s_1), n\equiv
s_1\bmod 4\mbox{ if } n< m_i(s_i)\}.$$
Observe that since $m_1(s_1)\leq m_2(s_1)$, $S_2^1(s_1)\subseteq S_2^2(s_1)$.

We now turn to the description of the third coordinate $s_3$. There are
three cases: $i=1, 1^*$ and $2$. We should also point out that in each
case, our description will break naturally into two parts. Borrowing
terminology from Wyman \cite{wyman}, there is {\em stable
ramification} when $s_3$ is uniquely determined by $s_1$ and $s_2$,
and there is {\em unstable ramification} when it is not.

We begin by describing $s_3$ under unstable ramification. In each
case, there are lower and upper bounds
$$L_i=\begin{cases}
7s_1-2s_2 &\mbox{for }i=1^*,\\
5s_1 &\mbox{for }i=1,\\
2s_1+3s_2 &\mbox{for }i=2,
\end{cases}\qquad
U_i=\begin{cases}
8e-3s_1 &\mbox{for }i=1^*, 1,\\
8e-2s_1-s_2&\mbox{for }i=2.
\end{cases}
$$ Notice that $L_{1^*}\leq L_1\leq L_2$ and $U_2\leq U_{1^*}=U_1$
(with equality everywhere, if we formally equate $s_1=s_2$).  When
there is room between the lower and upper bounds, namely $L_i<
U_i$, we have unstable ramification and
$$s_3\in\mbox{}^u\!S_3^i(s_1,s_2)=\{n\in \bZ: L_i\leq n\leq
U_i, s_3\equiv s_i\bmod 8\mbox{ if } L_i< n< U_i\}.$$
Note that in the description of $^u\!S_3^i(s_1,s_2)$, ``$s_3\equiv
s_i\bmod 8$'' means $s_3\equiv s_2\bmod 8$ for $i=2$ and $s_3\equiv
s_1\bmod 8$ for $i=1,1^*$.  Note furthermore that the condition
$L_i< U_i$, means $5s_1-s_2< 4e$ for $i=1^*$, $s_1< e$ for
$i=1$ and $s_1+s_2< 2e$ for $i=2$. Outside of this condition we
have stable ramification:
$$s_3=4e+\begin{cases}
s_i&\mbox{for }i=1,2,\\
2s_1-s_2&\mbox{for }i=1^*,
\end{cases}$$
which, of course, defines a set $^s\!S_3^i(s_1,s_2)$ for each
$i\in\{1,1^*,2\}$.

In summary, we have defined three sets, $\mathcal{R}_{1^*}^e,
\mathcal{R}_1^e,
\mathcal{R}_2^e$:
$$\mathcal{R}_i^e=\left\{(s_1,s_2,s_3)\in\bZ^3: s_1\in S_1, s_2\in
S_i(s_1), s_3\in \mbox{}^u\!S_3^i(s_1,s_2)\cup
\mbox{}^s\!S_3^i(s_1,s_2)\right\}.$$

Our interest in these sets is largely a result of our interest in
counter-examples to the conclusion of Hasse-Arf. By a result of
Fontaine \cite[{Prop 4.5}]{fontaine}, we know that such
counter-examples can occur only in the one break case. So we focus now
on $\mathcal{R}_{1^*}^e, \mathcal{R}_1^e$. We are principally
interested in the relationship between $s_1$ and $s_3$, which as we
might suspect from \S1.2 corresponds to the two (usual) ramification
breaks in a quaternion extension. To provide a two dimensional visual
aid, we slice now each of $\mathcal{R}_1^e, \mathcal{R}_{1^*}^e$ by
the hyperplane $s_2-s_1=e/2$, and project each slice to the
$(s_1,s_3)$-plane (with axes scaled 1--2).  The result (sketched
below) includes a line segment (representing stable ramification)
along with a triangular region (representing unstable ramification).
To aid comparison, we have included certain dotted segments of the
lines $s_3=3s_1$, $s_3=5s_1$, $s_3=s_1+4e$ and $s_3=8e-3s_1$ in both
sketches.  Note that since the upper bound for $s_2$, namely
$\min\{2s_1,4e-s_1\}$ depends upon whether or not $s_1\leq 4e/3$, the
hyperplane $s_2-s_1=e/2$ intersects $\mathcal{R}_1, \mathcal{R}_1^*$
only for $e/2\leq s_1 \leq 7e/4$.

\vspace*{3mm}
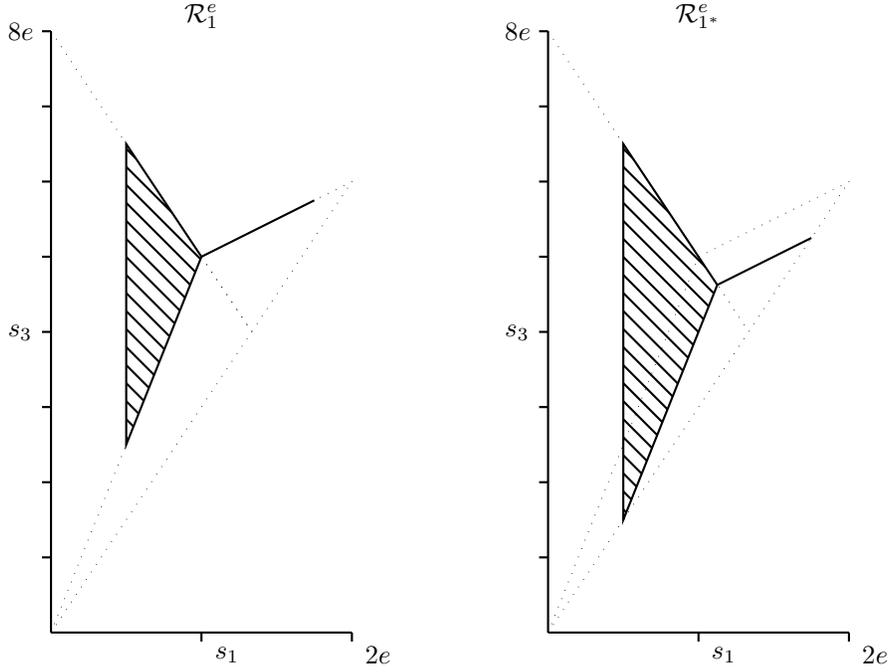
\begin{center}
\begin{tabular}{ccc}
$\mathcal{R}_1^e$ & & $\mathcal{R}_{1^*}^e$\\
\begin{pspicture}(0,0)(4,8)
\psline[linestyle=solid](2,5)(3.5,5.75)
\psline[linestyle=solid](0,0)(0,8)
\psline[linestyle=dotted, linewidth=.4pt](2,5)(4,6)
\psline[linestyle=dotted, linewidth=.4pt](2,5)(2.67,4)
\psline[linestyle=dotted, linewidth=.4pt](0,8)(2.67,4)
\psline[linestyle=dotted, linewidth=.4pt](2,5)
\psline[linestyle=dotted, linewidth=.4pt](4,6)
\pspolygon[fillstyle=vlines](1,2.5)(2,5)(1,6.5)
\psaxes[axesstyle=frame,labels=none,Dx=2,tickstyle=bottom](0,0)(0,8)(4,0)
\uput[0](-.75,4){$s_3$}
\uput[0](-.75,8){$8e$}
\uput[0](2,-.3){$s_1$}
\uput[0](4,-.3){$2e$}
\end{pspicture}
&\hspace*{.75in}&
\begin{pspicture}(0,0)(4,8)
\psline[linestyle=solid](2.25,4.625)(3.5,5.25)
\psline[linestyle=solid](0,0)(0,8)
\psline[linestyle=dotted, linewidth=.4pt](2,5)
\psline[linestyle=dotted, linewidth=.4pt](0,8)(2.67,4)
\psline[linestyle=dotted, linewidth=.4pt](2,5)(4,6)
\psline[linestyle=dotted, linewidth=.4pt](4,6)
\pspolygon[fillstyle=vlines](1,1.5)(2.25,4.625)(1,6.5)
\psaxes[axesstyle=frame,labels=none,Dx=2,tickstyle=bottom](0,0)(0,8)(4,0)
\uput[0](-.75,4){$s_3$}
\uput[0](-.75,8){$8e$}
\uput[0](2,-.3){$s_1$}
\uput[0](4,-.3){$2e$}
\end{pspicture}
\end{tabular}
\end{center}
\vspace*{6mm}

\subsection{Statement of Main Results}
Given a base field $K$, there is a map that sends each fully ramified
quaternion extension $N/K$ to its ramification triple, either
$(b_1,b_2,b_3)$ or $(b,r,b_3)\in \bZ^3$ depending upon the filtration
in its biquadratic subfield $M/K$.  A catalog of such ramification
triples should be considered complete if
\begin{enumerate}
\item the range of this map is given an explicit
description, and
\item the map is also shown to be onto this range.
\end{enumerate}
Using this definition, the catalog of ramification triples that we
give below is complete for fields $K$ that contain the $4$th roots of
unity, namely $\sqrt{-1}\in K$ ({\em i.e.} $K$ has level one).
Otherwise, when $\sqrt{-1}\not\in K$ our results are not
complete\footnote{Assumptions on the roots of unity in the base field
are common in ramification theory.  For example, the ramification
break number $b$ of a ramified $C_p$-extension $L/K$ of local number
fields with residue characteristic $p$, satisfies $1\leq b\leq
pe_K/(p-1)$ with the additional condition that $\gcd(b,p)=1$ for
$b<pe_K/(p-1)$.  But the case $b=pe_K/(p-1)$ is possible only when $K$
contains a $p$th root of unity \cite[III\S2 Prop 2.3]{fesenko}. As a
further example, note that \cite[Thm 32]{wyman} concerning
$C_{p^2}$-extensions is proven under the assumption that the base
field contains the $p$th roots of unity.} -- in particular, we do not
attempt to address condition (2).

In \S1.2, we decomposed the collection of fully ramified quaternion
extensions of $K$ into three subclasses: $\mathcal{Q}_{1^*}^K,
\mathcal{Q}_1^K, \mathcal{Q}_2^K$.  In \S1.3, we defined three ranges
(subsets of $\bZ^3$): $\mathcal{R}_{1^*}^{e}, \mathcal{R}_1^{e},
\mathcal{R}_2^{e}$.

\begin{theorem}
If $\sqrt{-1}\in K$ and $N/K$ is a fully ramified quaternion
extension, so $N\in \mathcal{Q}_i^K$ for some $i\in\{1,1^*,2\}$. Then
its ramification triple, either $(b,r,b_3)$ or $(b_1,b_2,b_3)$, lies
in $\mathcal{R}_i^{e_K}$, where $e_K$ denotes the absolute
ramification of $K$.

Moreover, given $K/\bQ_2$ with $\sqrt{-1}\in K$ and any
triple $(s_1,s_2,s_3)\in\mathcal{R}_i^{e_K}$ where
$i\in\{1,1^*,2\}$, there is a fully ramified quaternion extension
$N/K$ with $N\in \mathcal{Q}_i^K$ whose ramification triple is
$(s_1,s_2,s_3)$.
\end{theorem}

\begin{theorem}
If $\sqrt{-1}\not\in K$ and $N/K$ is a fully ramified quaternion
extension with $N\in \mathcal{Q}_1^K\cup \mathcal{Q}_2^K$ and a stable
ramification triple, either $(b,r,b_3)$ with $b>e_K$ or
$(b_1,b_2,b_3)$ with $b_1+b_2>2e_K$, then its ramification triple lies
in $\mathcal{R}_i^{e_K}$, where $e_K$ denotes the absolute
ramification of $K$.
\end{theorem}

\subsection{Hasse-Arf} The Hasse-Arf Theorem
states that upper ramification numbers in abelian extensions are
integers.  Our results confirm a result of Fontaine, which says that
in quaternion extensions upper ramification numbers generally are
integers.

\begin{theorem}[{\cite[{Prop 4.5}]{fontaine}}]
Upper ramification numbers of fully ramified quaternion extensions of
dyadic number fields are integers, except when there is only one break
$b$ for $M/K$ and $b_3=3b$.
\end{theorem}
Note that the upper ramification numbers in a
quaternion extension are integers precisely when $b_3\equiv b_2 \mbox{
(or }b) \bmod 4$ \cite[IV\S3]{serre:local}.
Using Theorems 1.1 and the illustration for
$\mathcal{R}_{1^*}^e$, we can see that the exceptional situation
$b_3=3b$ {\em can occur} only when both the second refined break is
{\em maximal}: $r=\min\{4e_K-b,2b\}$, and its associated root of unity
$\omega$ is {\em minimal}: $\omega^3=1$.  Based upon
\S1.4, $b_3=3b$ {\em must occur} under $r=4e_K-b$ and
$\omega^3=1$.

\section{Embeddability and Quadratic Defects}

In 1936 E.~Witt characterized the biquadratic extensions
$M=K(\sqrt{u},\sqrt{v})$ that embed in a quaternion extension
\cite{witt}.  When $K$ is a finite extension of $\bQ_2$, his condition
is equivalent to the Hilbert symbol equality: $(-u,-v)=(-1,-1)$, which
is equivalent to the product formula $(-1,u)(-1,v)(u,v)=1$.

If the product formula holds then, replacing
$u$ or $v$ with $uv$ if necessary and using Hilbert symbol properties,
we may assume without loss of generality that $(u,v)=1$ and
$(uv,-1)=1$. As a result, when $M$ embeds in a quaternion extension,
we may assume that there are two elements $\eta\in K(\sqrt{u})$ and
$\tau\in K(\sqrt{uv})$ whose norms satisfy $N_{K(\sqrt{u})/K}(\eta)=v$
and $N_{K(\sqrt{uv})/K}(\tau)=-1$.

An observation of H.~Reichardt then characterizes the quaternion
extensions $N/K$ that contain $M$: for if we let
$\alpha_k\in M$ be defined by
$$\alpha_k=k\cdot \sqrt{uv}\cdot\eta\cdot\begin{cases}1
&i=\sqrt{-1}\in K\\ \tau &i=\sqrt{-1}\not\in K\end{cases}$$ where
$k\in K$ and $N=M(\sqrt{\alpha_k})$, then $N/K$
is a quaternion extension \cite{reichardt}.  Moreover it is generic in
the sense that any quaternion extension of $K$ containing $M$ can be
expressed as $M(\sqrt{\alpha_k})$ for some $k\in K$.  Jensen and Yui
provide a nice source for these results. Indeed Witt's condition
\cite[Lem I.1.1]{jensen:yui} can be translated to the Hilbert symbol
condition using \cite[Lem I.1.6]{jensen:yui}, and Reichardt's
observation appears as \cite[Lem I.1.2]{jensen:yui}.

\subsection{Quadratic Extensions and Quadratic Defect}

Let $F$ denote a finite extension of $\bQ_2$, and let $T_F$ be its
maximal unramified subfield.  As is well-known, a vector space basis
for $F^*/(F^*)^2$ over $\bF_2$ is given by $\{1+a\pi_F^{2n-1}:
a\in\euO_{T_F}/2\euO_{T_F}, 1\leq n\leq e_F\}$ along with $\pi_F$ and
$1+4\lambda$ for some $\lambda\in \euO_{T_F}$ with $x^2+x+\lambda$
irreducible over ${T_F}$, \cite[Ch15]{hasse}. It is easy to check that
$F(\sqrt{1+4\lambda})/F$ is unramified. This means that there are
essentially two types of ramified quadratic extensions: those that
arise from the square root of a prime, $F(\sqrt{\pi_F})$, and those
that arise from the square root of a one-unit, $F(\sqrt{u})$ with
$u=1+\beta$ and $0<v_F(\beta)<2e_F$ odd. Define the defect in $F$ of a
prime element to be $\mbox{def}_F(\pi_F)=0$ and of a unit to be
$\mbox{def}_F(u)=\max\{v_F(k^2u-1):k\in F\}$ \cite[\S63A]{omeara}.  If
$u=1+\beta$ as above, $\mbox{def}_M(1+\beta)=v_F(\beta)$.  It is
straightforward now to verify that the ramification number of
$F(\sqrt{\kappa})/F$ (for $\kappa\in F^*\setminus (F^*)^2$) is tied to
the defect of $\kappa$ by $b=2e_F-\mbox{def}_F(\kappa)$. (All this is
generalized to include odd primes $p$ in \cite[\S4]{wyman}.)

Recall Question 1.  Given a quaternion extension $N/K$, we are
interested in determining $b_3$, the ramification break for the
quadratic extension $N/M$, which is tied to the quadratic defect of
$\alpha_k$ in $M$ by
$$b_3=8e_K-\mbox{def}_M(\alpha_k).$$ Indeed we will determine $b_3$ by
determining $\mbox{def}_M(\alpha_k)$.  Recall that $\alpha_k$ is a
product: either $k\cdot \sqrt{uv}\cdot\eta$ or $k\cdot
\sqrt{uv}\cdot\eta\cdot\tau$ depending upon whether $\sqrt{-1}\in K$
or not.  It is easy to see that $\mbox{def}_M(A\cdot
B)\geq\min\{\mbox{def}_M(A),\mbox{def}_M(B)\}$ for $A, B\in M$, and
that we can be certain of equality only when $\mbox{def}_M(A) \neq
\mbox{def}_M(B)$.  The technical work in this paper addresses two
issues: (1) The terms in $\alpha_k$ lie in proper subfields of $M$. As
a result, the defect in $M$ of each term is not immediately obvious
from its expression.  (2) Moreover, once the defect of each term has
been determined, there are often at least two terms with the same
defect.

\subsection{Two Technical Lemmas}
If $E/F$ is a ramified quadratic extension and $\kappa\in F$, then
$\mbox{def}_E(\kappa)> \mbox{def}_F(\kappa)$. To describe this
increase in valuation carefully, we need to define the following
continuous increasing function.
$$g_{F,b}(x)=\min\{ 2x+b,
x+ 2e_F\}=\begin{cases}2x+b&\mbox{for } x\leq 2e_F-b,\\
x+2e_F&\mbox{for } x> 2e_F-b.\end{cases}$$

\begin{lemma}
Let $E/F$ be a ramified quadratic extension with break number $b$ odd.
If $\kappa\in F^*\setminus(E^*)^2$, then
$$\mbox{\rm def}_E(\kappa)\geq g_{F,b}\left (\mbox{\rm
def}_F(\kappa)\right )$$ with equality when $\mbox{\rm
def}_F(\kappa)\neq 2e_F-b$. As a result, given a threshold value
$\delta\geq 0$ with $\mbox{\rm def}_F(\kappa)\geq \delta$, then
$\mbox{\rm def}_E(\kappa)\geq g_{F,b}(\delta) \geq b$.
\end{lemma}
\begin{proof}
Since $\kappa\not\in(E^*)^2$, $E(\sqrt{\kappa})/F$ is a biquadratic
extension. The break number of $F(\sqrt{\kappa})/F$ is
$2e_F-\mbox{def}_F(\kappa)$. Passing to the upper numbering
for the filtration of $E(\sqrt{\kappa})/F$
\cite[ChIV\S3]{serre:local}, we see that the break number for
$E(\sqrt{\kappa})/E$ is
$$
\begin{array}{cl}
 4e_F-2\mbox{def}_F(\kappa)-b & \quad \textrm{for \ } 2e_F-\mbox{def}_F(\kappa)>b\,; \\
 2e_F-\mbox{def}_F(\kappa) & \quad \textrm{for \ } 2e_F-\mbox{def}_F(\kappa)<b \,; \\
 \leq b & \quad \textrm{for \ } 2e_F-\mbox{def}_F(\kappa)=b\,.
\end{array}
$$
As a result,
$$
\mbox{def}_E(\kappa) = \left\{
\begin{array}{cl}
 2\mbox{def}_F(\kappa)+b & \quad \textrm{for \ } 2e_F-\mbox{def}_F(\kappa)>b\,; \\
 2e_F+\mbox{def}_F(\kappa) & \quad \textrm{for \ } 2e_F-\mbox{def}_F(\kappa)<b \,; \\
 \geq 4e_F-b & \quad \textrm{for \ } 2e_F-\mbox{def}_F(\kappa)=b\,.
\end{array} \right.
$$
The result follows. \end{proof}

Given elements of known defect, can one be described in terms of
another?
\begin{lemma}
Given $\beta\in F$ with $v_F(\beta)=2e_F-b$ and $0<b<2e_F$ odd.  If
$\kappa\in 1+\euP_F$ and $\mbox{\rm def}_F(\kappa)=2e_F-a$ with
$0<a\leq b$, there is a $\mu\in\euO_F$ with $v_F(\mu)=(b-a)/2$ and a
$\lambda\in \euO_{T_F}$, either $0$ or so that $z^2+z=\lambda$ is
irreducible over ${T_F}$, such that
$$\kappa=(1+\mu^2\beta)(1+4\lambda)\in(1+\euP_F)/(1+\euP_F)^2.$$
\end{lemma}
\begin{proof}

Clearly $\kappa=
(1+\mu_0^2\beta)(1+\mu_1^2\beta)(1+\mu_2^2\beta)\cdots (1+4\lambda)
\in(1+\euP_F)/(1+\euP_F)^2$ for some $\mu_i\in \euO_F$ with
$(b-a)/2=v_F(\mu_0)<v_F(\mu_1)<v_F(\mu_2)<\cdots$. Now use
$(1+\mu_0^2\beta)(1+\mu_1^2\beta)\equiv 1+(\mu_0+\mu_1)^2\beta\bmod
\mu_1^2\beta\euP_F$ repeatedly.
\end{proof}

\section{One Break Biquadratic Extensions}
Let $M/K$ be a fully ramified biquadratic extension which has only one
ramification break, at $b$. In this case the ramification numbers for
each of the three subfields must be the same. Using Lemma 2.2, there
must be a $\beta\in K$ with $v_K(\beta)=2e_K-b$; a nontrivial
$2^f-1=q-1$ root of unity $\omega\in\euO_T$; a $\mu\in K$ where either
$\mu=0$ or $v_K(\mu)=m$ with $0< m<b/2$; and a $\lambda\in K$ where
either $\lambda=0$ or $\lambda$ is a $q-1$ root of unity with
$z^2+z=\lambda$ irreducible over $T$; such that $M=K(x,y)$ where
\begin{eqnarray*}
x^2&=&1+\beta,\\
y^2&=&\left (1+(\omega+\mu)^2\beta\right )(1+4\lambda).
\end{eqnarray*}
Without loss of generality, we let $\sqrt{u}=x$ and $\sqrt{v}=y$, and
for the remainder of this section set $L=K(x)$. Note that because
$N_{L/K}(x-1)=-\beta$ and $v_K(\beta)=2e_K-b$, we must have
$v_L(x-1)=2e_K-b$.  We let
$\overline{G}=\mbox{Gal}(M/K)=\langle\sigma, \gamma\rangle$ where the
generators act by $\sigma x=x $ and $\gamma y=y$. (It should cause no
confusion that we use $\sigma, \gamma$ to denote both the generators
of the quaternion group $G$ and its $C_2\times C_2$ quotient group
$\overline{G}$.)

Before we turn to the refined ramification filtration, notice that the
extension $M/L$ is quadratic with break $b$. As a result, there should
be a unit $U\in M$ of defect $\mbox{def}_M(U)=4e_K-b$ such that
$M=L(U)$. Motivated by and identity in $\mathbb{Q}(A,X)$,
\begin{equation}
(1+A(X-1))^2=(1+A^2B) \left(1+2(A-A^2)(X-1)\frac{1}{1+A^2B}\right ),
\end{equation}
where $B=X^2-1$, we choose $Y\in M$ so that
\begin{equation}
yY=1+(\omega+\mu)(x-1)\in L.
\end{equation}
Now using (1) with $X=x$, $B=\beta$, and
$A=\omega+\mu$, we find that
\begin{equation}
Y^2=\left(1+2((\omega+\mu)-(\omega+\mu)^2)(x-1)
\frac{1}{1+(\omega+\mu)^2\beta}\right ) (1+4\lambda)^{-1}.
\end{equation}
As a result, by applying the norm $N_{M/L}(Y-1)=1-Y^2$
where $v_L(1-Y^2)=4e_K-b$. Thus
$v_M(Y-1)=4e_K-b$. So $Y$ is our desired
unit, and $M=L(Y)$ with $\sigma Y=-Y$.

\subsection{The second refined break and its associated root of unity}

When there is only one ramification break, all Galois action
``looks'' the same from the perspective of the usual ramification
filtration. Thus the necessity of a {\em refined ramification
filtration}, which helps us ``see'' a difference. As an aid to the
reader, we replicate some of the material from
\cite{elder:newbreaks}, restricting to $p=2$, so that many of the
details are simpler.

Let $J=(\sigma-1, \gamma-1)$ be the Jacobson
radical of $\bF_q[\overline{G}]$.
Define an $\bF_q$-`action' on the
one-units $1+J$ by the map
$$(a,1+x)\in\bF_q\times (1+J)\longrightarrow x^{[a]}:=1+ax\in 1+J.$$
This makes $1+J$ a near space over $\bF_q$ with all the properties of
a vector space, except that scalar multiplication does not necessarily
distribute: It is possible to find $x,y\in J$ and $a\in \bF_q$ so that
$((1+x)(1+y))^{[a]}(1+x)^{[-a]}(1+y)^{[-a]}=1+(a^2+a)xy\neq 1$.
We do not have a proper action. To achieve one and
create a vector space, we deviate slightly from \cite{elder:newbreaks}
and define $$\overline{G}^{\mathcal{F}}:=(1+J)/(1+J^2).$$ It is
straightforward to check that this vector space over $\bF_q$ has
basis $\{\sigma, \gamma\}$.

To define a ramification filtration for $\overline{G}^{\mathcal{F}}$,
choose any element $\rho\in M$ with valuation $v_M(\rho)=b$, and
define, for $s\in \overline{G}^{\mathcal{F}}$,
$w_{\rho}(s)=\max\{v_M((\tilde{x}-1)\rho):\tilde{x}\in\euO_T[\overline{G}],
x=\tilde{x}+2\euO_T[\overline{G}], x\in\bF_q[\overline{G}],
s=x(1+J^2)\}$, and the refined ramification groups by
$$\overline{G}^{\mathcal{F}}_i=\{s\in \overline{G}^{\mathcal{F}}:w_{\rho}(s)-v_M(\rho)\geq
i\}.$$

For example, we will use $\rho=2/(Y-1)$. If we replace $y^2$ with
$(\omega^{-1}y)^2=(\omega^{-2}+\beta)(1 + \mu^2\beta)\bmod
\mu^2\beta\euP_K$, then we have notation (including $\sigma, \gamma$ as
generators of the Galois group) exactly as in
\cite[\S4.1]{elder:newbreaks}. As a result, we can apply
\cite[Prop4.2]{elder:newbreaks} and find that
$$v_M((\gamma\sigma^{[\omega]}-1)\rho)=v_M(\rho)+r$$ where
$r=\min\{4e_K-b,b+4m, 2b\}$. Thus there are two breaks in the refined
filtration: namely $b<r$ with $\overline{G}^{\mathcal{F}}_b\supsetneq
\overline{G}^{\mathcal{F}}_{b+1}$ and
$\overline{G}^{\mathcal{F}}_r\supsetneq
\overline{G}^{\mathcal{F}}_{r+1}$, where the {\em second refined
break} satisfies $r\leq \min\{2b,4e_k-b\}$ and away from this upper
bound satisfies $r\equiv b\bmod 4$.  These breaks are independent of
our choices: of $\rho$ and of the generators for $\overline{G}$.

Additionally, the second refined break $r$ is associated with a root of
unity, namely $\omega$, which does depends upon our choice of
generators for $\overline{G}$.  Replace $\gamma$ by $\gamma\sigma$ and
we have an alternative root of unity $\equiv\omega+1\bmod 2$. Indeed
these are the only two roots of unity that arise from a change of
generators for $\mbox{Gal}(M/K)$. This suggests an equivalence
relation on nontrivial $q-1$ roots of unity: $\omega\sim\omega'$ if
and only if $\omega\equiv\omega'$ or $\omega'+1\bmod 2$.  If we
identify these nontrivial $q-1$ roots of unity with their images in
$\bF_q\setminus\bF_2$, then the equivalence classes of this relation
can be identified with the $f_K-1$ nontrivial additive cosets of
$\bF_q/\bF_2$.  Thus the second refined break $r$ is actually
associated with an equivalence class of two $q-1$ roots of unity. We
are going to be interested in whether the elements of a particular
equivalence class satisfy a condition: whether they both are
nontrivial cube roots of unity. So it is worth pointing out that
$\omega^2+\omega+1=0\bmod 2$ if and only if
$(\omega+1)^2+(\omega+1)+1\equiv 0\bmod 2$. As a result, in the
statements of our results we can refer to ``$r$ and its associated
root of unity $\omega$'' (equating each equivalence class with a
representative). The condition $\omega^3=1$ is well-defined.

\begin{remark}
In the two break case, this refined ramification filtration produces
the usual two ramification break numbers.
\end{remark}

\subsection{The determination of $\mbox{def}_M(\alpha_k)$}
We begin with a lemma.
\begin{lemma}
Let $M/K$ be a fully ramified biquadratic extension with one
ramification break, at $b$.  Then for all $k\in K$, $\mbox{\rm
def}_M(k)\geq 3b$.  In particular, $\mbox{\rm def}_M(\pi_K)=3b$.  And
if $k\in 1+\euP_K$ with $0<\mbox{\rm def}_K(k)<2e_K-b$, then
$\mbox{\rm def}_M(k)=3b+4\mbox{\rm def}_K(k)$. So $3b<\mbox{\rm
def}_M(k)<8e_K-b$ with $\mbox{\rm def}_M(k)\equiv -b\bmod 8$.
Otherwise if $k\in 1+\euP_K$ with $2e_K-b\leq \mbox{\rm def}_K(k)$
then $\mbox{\rm def}_M(k)\geq 8e_K-b$.
\end{lemma}
\begin{proof}
Apply Lemma 2.1. Note that for $0<\mbox{\rm def}_K(k)<2e_K-b$,
$\mbox{\rm def}_K(k)\equiv b \bmod 2$, therefore $3b+4\mbox{\rm
def}_K(k)\equiv 7b\bmod 8$.
\end{proof}

Now we assume that $M=K(x,y)$ embeds in a quaternion extension, that
$x=\sqrt{u}, y=\sqrt{v}$, and adopting notation as in \S2, we
determine $\mbox{def}_M(\alpha_k)$. Our analysis decomposes: Case $1$ in
\S3.3 when $\omega^3\neq 1$.  and Case $1^*$ in \S3.4 when $\omega^3= 1$.
We begin with the easier case.

\subsection{\textbf{Case $1$: Assume $\omega^3\neq 1$.}}
Under this assumption there is both stable and unstable
ramification. We begin with stable situation.

\subsubsection{\textbf{Stable ramification: $b>e_K$}}
We {\em do not} assume $\sqrt{-1}\in K$.
It follows immediately from the following lemma that if $b>e_K$ then
$b_3=4e_K+b$, the value given in the catalog in
Section~\ref{ssec:catalog}.
\begin{lemma}
If $b>e_K$ and $\omega^3\neq
1$, then $\mbox{\rm def}_M(\alpha_k)=4e_K-b$.
\end{lemma}
\begin{proof}
Because of the possibility that $i=\sqrt{-1}\not\in K$, we have
$\alpha_k=kxy\eta\tau$.  It suffices therefore to check that
$\mbox{def}_M(k\cdot yY\eta\cdot \tau)>4e_K-b$ and that
$\mbox{def}_M(xY)=4e_K-b$.

We prove the first statement, $\mbox{def}_M(k\cdot yY\eta\cdot
\tau)>4e_K-b$, by showing $\mbox{def}_M(k)>4e_K-b$,
$\mbox{def}_M(yY\eta)>4e_K-b$ and $\mbox{def}_M(\tau)>4e_K-b$.  Using
Lemma 3.1 and $b>e_K$, we have $\mbox{def}_M(k)\geq 3b>4e_K-b$. The
other two inequalities will follow from Lemma 2.1 if we can show
$\mbox{def}_L(yY\eta)>2e_K-b$ and $\mbox{def}_{K(xy)}(\tau)>2e_K-b$,
respectively.  Recall that $L=K(x)$. We begin with
$\mbox{def}_L(yY\eta)$: Note that by (2), $yY\in L$ and by assuming
$x=\sqrt{u}$, $\eta\in L$ as well.  Now check, using $b>e_K$, that
$N_{L/K}(1+\omega(x-1))\equiv 1+\omega^2\beta\bmod \beta\euP_K$.  This
implies that $\eta\equiv 1+\omega(x-1)\bmod (x-1)\euP_L$.  Hence using
(2) we see that $\mbox{def}_L(yY\eta)>2e_K-b$.  Finally consider
$\mbox{def}_{K(xy)}(\tau)$: Using \cite[V\S3]{serre:local} we see that
if $\tau\not\in (K(xy)^*)^2$, because $\mbox{def}_K(-1)\geq e_K$ we
have $\mbox{def}_{K(xy)}(\tau)\geq e_K$.  Hence, since $b>e_K$ we have
$\mbox{def}_{K(xy)}(\tau)>2e_K-b$.

We prove the second statement, $\mbox{def}_M(xY)=4e_K-b$.
Notice that $\mbox{def}_M(Y)=v_M(Y-1)=4e_K-b$ and that Lemma 2.1 gives
$\mbox{def}_M(x)=4e_K-b$. This makes $\mbox{def}_M(xY)$ difficult to
determine, but also means that there is a unit $a\in \euO_K$ such that
$x=1+a(Y-1)$ in
$\mathcal{M}:=(1+\euP_M)/[(1+\euP_M)^2(1+(Y-1)\euP_M)]$.
We will have the desired conclusion if we can show $a\not\equiv
1\bmod\euP_K$.  Since $v_L(2)=2e_K$ is even, there is a $\kappa\in L$
such that $\kappa^2\equiv 2(\omega-\omega^2)\bmod 2\euP_L$.  Using
(1), we expand $(1+(1/\kappa)(Y-1) )^2=A\cdot B\bmod (Y-1)\euP_M$ where
$A=1+(Y^2-1)/\kappa^2\in L$ and $B=1+2(1/\kappa-1/\kappa^2)(Y-1)$.
This means that $A=B$ in $\mathcal{M}$.  Since $(Y^2-1)/\kappa^2\equiv
(x-1)\bmod (x-1)\euP_L$, we see that $A=x\cdot C$ where
$\mbox{def}_L(C)>2e_K-b$ and thus by Lemma 2.1,
$\mbox{def}_M(C)>4e_K-b$.  So $A=x$ in $\mathcal{M}$, and thus $x=B$
in $\mathcal{M}$.  Notice that $B\equiv 1-2(Y-1)/\kappa^2\equiv
1+(\omega-\omega^2)^{-1}(Y-1)\bmod (Y-1)\euP_M$. This means that
$x=1+(\omega^2+\omega)^{-1}(Y-1)$ in $\mathcal{M}$. And because $\omega^3\neq 1$,
$(\omega^2+\omega)^{-1}\not\equiv 1\bmod\euP_K$.
\end{proof}

\subsubsection{\textbf{Two preliminary results:}}

For the remaining cases, we need two additional technical results.
Define the following:
\newline For $b<e_K$, (for use in \S3.3.3)
\begin{eqnarray*}
\mathcal{L}&=&\frac{1+\euP_L}{(1+\euP_K)(1+\euP_L)^2(1+(\beta/2)(x-1)\euP_L)},\\
\mathcal{M}&=&\frac{1+\euP_M}{(1+\euP_K)(1+\euP_M)^2(1+(\beta/2)(Y-1)\euP_M)}.
\end{eqnarray*}
For all $b$, (for use in \S3.4)
\begin{eqnarray*}
\mathcal{L}^*&=&\frac{1+\euP_L}{(1+\euP_K)(1+\euP_L)^2(1+\beta\euP_L)},\\
\mathcal{M}^*&=&\frac{1+\euP_M}{(1+\euP_K)(1+\euP_M)^2(1+(x-1)(Y-1)\euP_M)}.
\end{eqnarray*}
It is easy to see that, under $b<e_K$, the natural maps,
$\mathcal{L}^*\rightarrow\mathcal{L}$ and
$\mathcal{M}^*\rightarrow\mathcal{M}$, are surjective. Moreover, we
can define defects with respect to these groups in the natural way.
For example, for $\mu\in 1+\euP_M$,
$\mbox{def}_{\mathcal{M}^*}(\mu)=\max\{v_M(m-1):m=\mu\in
\mathcal{M}^*\}$. Note: We will regularly abuse notation by
identifying a coset with one of its coset representatives.

\begin{lemma}
The inclusions $1+(\beta/2)(x-1)\euP_L\subseteq
(1+\euP_M)^2(1+(\beta/2)(Y-1)\euP_M)$ for $b<e_K$,
$1+\beta\euP_L\subseteq (1+\euP_M)^2(1+(x-1)(Y-1)\euP_M)$ yield the
following natural, well-defined maps:
$\mathcal{L}\rightarrow\mathcal{M}$ defined for $b<e_K$ and
$\mathcal{L}^*\rightarrow\mathcal{M}^*$ defined for all $b$.
\end{lemma}
\begin{proof}
Use Lemma 2.1 to determine the two inclusions.
\end{proof}

\begin{lemma}
For $i=\sqrt{-1}\in K$, the coset equality holds for all
$b$:
$$Y= (1+i(\omega+\omega^2+\mu+\mu^2)(x-1))\cdot
(1+(\omega+\omega^2)(x-1)(Y-1))\in\mathcal{M^*}.$$ Moreover for
$b<e_K$, $Y= (1+i(\omega+\omega^2+\mu+\mu^2)(x-1)) \in\mathcal{M}$.
For $b>e_K$, $Y= (1+(\omega+\omega^2+\mu+\mu^2)(x-1))\cdot
(1+(\omega+\omega^2)(x-1)(Y-1))\in\mathcal{M^*}$.
\end{lemma}
\begin{proof}
Expand $(1+(Y-1)/(i-1))^2=1-(1+i)(Y-1)+i(Y-1)^2/2=
1-(1+2i)(Y-1)+i(Y^2-1)/2\equiv Y+i(Y^2-1)/2\bmod (x-1)(Y-1)\euP_M$,
noting that $v_M(2(Y-1))>v_M((x-1)(Y-1))$.  Factor
$Y+i(Y^2-1)/2=A\cdot B$ with
$$A=\left (1+i\frac{Y^2-1}{2}\right )\in L,\qquad B= \left (
1+(Y-1)\frac{1}{1+i\frac{Y^2-1}{2}}\right )\in M.$$ Each of $A$ and
$B$ has a copy of $i(Y^2-1)/2$ that needs to be replaced. From (3) we
have the approximation $(Y^2-1)/2\equiv
(\omega+\omega^2+\mu+\mu^2)(x-1)-2\lambda \bmod \beta\euP_L$. So
$A\equiv (1+i(\omega+\omega^2+\mu+\mu^2)(x-1))\cdot (1-2i\lambda)\bmod
\beta\euP_L$.
Because $b>e_K$ means $v_L((i-1)(x-1))> v_L(\beta)$,
we can drop the first ``$i$'' in this expression
when $b>e_K$.
Thus we find that as elements of $\mathcal{L}^*$, and
using Lemma 3.3, also as elements of $\mathcal{M}^*$,
$$A=\begin{cases}
1+(\omega+\omega^2+\mu+\mu^2)(x-1)&\mbox{for } b>e_K,\\
1+i(\omega+\omega^2+\mu+\mu^2)(x-1)&\mbox{for }b\leq e_K.
\end{cases}$$
We also have $i(Y^2-1)/2\equiv (\omega+\omega^2)(x-1)\bmod
(x-1)\euP_L$, which yields $B \equiv Y\cdot
(1+(\omega+\omega^2)(x-1)(Y-1))\bmod (x-1)(Y-1)\euP_M$.  So as
elements of $\mathcal{M}^*$, we also have $B= Y\cdot
(1+(\omega+\omega^2)(x-1)(Y-1))$. And by putting everything together,
we get the result.
\end{proof}

\subsubsection{\textbf{Unstable ramification: $b\leq e_K$}}
Assume that $i=\sqrt{-1}\in K$. Then
$e_K$ must be even. But since $b$ is odd, this means that we are really
assuming $b<e_K$.

In the following lemma we prove that if $i\in K$, $b< e_K$ and
$\omega^3\neq 1$, there is a $k_0\in K$ such that
$\mbox{def}_M(\alpha_{k_0})= 8e_K-5b$. From Lemma 3.1 it therefore
follows that $3b\leq \mbox{def}_M(\alpha_k)\leq 8e_K-5b$, and that
$\mbox{def}_M(\alpha_k)\equiv -b\bmod 8$ when $3b<
\mbox{def}_M(\alpha_k)< 8e_K-5b$. The values for $b_3$ listed in the
catalog in Section~\ref{ssec:catalog} follow immediately. Moreover,
each of these values for $\mbox{def}_M(\alpha_k)$ is realized.

\begin{lemma}If $i\in K$, $b< e_K$ and
$\omega^3\neq 1$, then
$\mbox{\rm def}_{\mathcal{M}}(xy\eta)= 8e_K-5b$.
\end{lemma}
\begin{proof}
Recall that because $i\in K$ we have $\alpha_k=kxy\eta$.  Since
$v_M((\beta/2)(Y-1))=8e_K-5b$, it is clear that our goal should be to
find a unit $u\in M$ such that $xy\eta = 1+u(\beta/2)(Y-1)$ as
elements in $\mathcal{M}$.  But since it is easier to work in $L$, we
first find an equivalent expression in $\mathcal{L}$ for $x\cdot
yY\cdot \eta\in L$. Then we use Lemma 3.4 to replace $Y$.

Note that because $b<e_K$, $v_K(\beta/2)=v_L((x-1)/(i-1)>0$.  Now
expand $(1+(X/(i+1))(x-1))^2$, using $\beta=x^2-1$, and find that for
$X\in \euO_K$,
\begin{equation}
1+ \frac{X^2(x-1)-(1+i)(X+X^2)(x-1)}{1+iX^2\beta/2}\in (1+\euP_K)(1+\euP_L)^2.
\end{equation}

If we
substitute $X=1$ in (4) and notice that
$v_L(2(1+i)(x-1))>v_L(\beta)$ and
$v_L((i+1)(\beta/2)(x-1))>v_L(\beta)$, we see that
$1+(x-1)/(1+\beta/2)\in (1+\euP_K)(1+\euP_L)^2(1+\beta\euP_L)$.
Therefore
\begin{equation}
x=1+(\beta/2)(x-1)\frac{1}{1+\beta/2}\mbox{\; in\; } \frac{1+\euP_L}{(1+\euP_K)(1+\euP_L)^2(1+\beta\euP_L)}.
\end{equation}
and since $v_L(\beta)>v_L((\beta/2)(x-1))$, we have
$x=1+(\beta/2)(x-1)\in
\mathcal{L}$.

Now substitute $X=\omega+\mu$ in (4). Simplify, again using
$v_L(2(1+i)(x-1))>v_L(\beta)=v_L((x-1)^2)>v_L((\beta/2)(x-1))$. This
results in the identity
$(1+(\omega+\mu)(x-1))\cdot(1-i(\omega+\mu+\omega^2+\mu^2)(x-1))\cdot(1+\omega^4(\beta/2)(x-1))=1$
in $\mathcal{L}$. Recall (2), namely
$yY=1+(\omega+\mu)(x-1)$. Therefore since $x=1+(\beta/2)(x-1)\in
\mathcal{L}$,
$$xyY=(1-i(\omega+\mu+\omega^2+\mu^2)(x-1))\cdot(1+(1+\omega^4)(\beta/2)(x-1))\in\mathcal{L}.$$

Using \cite[V\S3]{serre:local}, and the fact that
$v_K((\omega+\mu)^2\beta)>b$ we can choose $\eta$ so that
$\mbox{def}_L(\eta)=4e_0-3b$. Note
$\mbox{def}_L(1+(\omega+\mu)^2(\beta/2)(x-1))=4e_0-3b$ and moreover
that the norm $(1+(\omega+\mu)^2(\beta/2)(x-1))^{\sigma+1}\equiv
1+(\omega+\mu)^2\beta\bmod 2$.  Since $e_K>b$, every element of
$1+2\euO_K$ is a norm from $L=K(x)$ of an element of defect $\geq
4e_K-b>4e_K-3b$. Therefore
$\eta=1+(\omega+\mu)^2(\beta/2)(x-1)=1+\omega^2(\beta/2)(x-1) \in
\mathcal{L}$. So
$$xyY\eta=(1-i(\omega+\mu+\omega^2+\mu^2)(x-1))\cdot(1+(1+\omega^2+\omega^4)(\beta/2)(x-1))\in\mathcal{L}.$$

Now using Lemma 3.4, we have
$$xy\eta=(1+(1+\omega^2+\omega^4)(\beta/2)(x-1))\in\mathcal{M}.$$
Note that $1+(1+\omega^2+\omega^4)(\beta/2)(x-1)\in L$ and
because $\omega$ is not a third root of unity, by Lemma
2.1 we see that
$\mbox{def}_M(1+(1+\omega^2+\omega^4)(\beta/2)(x-1))=8e_K-5b$.
\end{proof}

\subsection{\textbf{Case $1^*$: Assume $\omega^3=1$.}}
Throughout this section we assume that $i=\sqrt{-1}\in K$.  Because of
$1+\omega+\omega^2=0$, we will require descriptions of
$\alpha_k=kxy\eta$ up to terms that have valuation strictly greater
than $8e_K-3b$. In other words, we will need to identify $\alpha_k$ in
$\mathcal{M}^*$.  This bound of $8e_K-3b=v_M((x-1)(Y-1))$ is
significantly larger than the bounds required in \S3.3.1 and \S3.3.3:
namely, $b+4e_K$ in the stable case and $8e_K-5b$ in the unstable
case. And this results in additional technicalities.

The material in \S3.3.3 is a good source of motivation.  Indeed, it
suggests that we proceed in two steps: First, identify an equivalent
expression in $\mathcal{L}^*$ for $x\cdot yY\cdot\eta\in L$. Most of
our technical difficulties are associated with the expression for
$\eta$.  Second, use Lemma 3.4 to replace $Y$.  There will be three
cases: (1) $b\leq e_K$, (2) $e_K<b<e_K+m$ and (3) $e_K+m\leq b$, each
associated with a different expression for $\eta$.  But in order to
keep the parallels to \S3.3.1 and \S3.3.3 evident, we present the
material in two sections: $b>e_K$, which is mostly stable
ramification, and $b\leq e_K$, which is most of unstable ramification.

\subsubsection{\textbf{Mostly stable ramification: $b>e_K$.}}
Using Lemma 3.1 and Lemma 3.7 below, we find that $\mbox{\rm
def}_M(kxy\eta)=\min\{4e_K-b+4m, 8e_K-3b\}$ for all $k\in K$ and
$b>e_K+m$. This is because $b>e_K+m$ can be rewritten as
$3b>4e_K-b+4m$ and so by Lemma 3.1, $\mbox{def}_M(k)\geq
3b>\mbox{\rm def}_M(k_0xy\eta)$ for all $k\in K$.  We also find that
for $e_K<b< e_K+m$, we have $3b\leq \mbox{\rm def}_M(kxy\eta)\leq
\min\{8e_K-5b+8m,8e_K-3b\}$ Moreover, each of these possible values
for $\mbox{def}_M(\alpha_k)$ is realized.

We start with a result that describes $\eta$.
\begin{lemma}
If $i=\sqrt{-1}\in K$, there is an $\eta^*\in L$ that satisfies
$$\mbox{\rm def}_L(\eta^*)
=\begin{cases}
2m+2e_K-b&\mbox{for }m+e_K<b,\\
4m+4e_K-3b&\mbox{for }m+e_K\geq b,
\end{cases}$$
and the
coset identity $\eta=(1+\omega(x-1))\cdot\eta^*$ in $\mathcal{L}^*$.
Furthermore $\mbox{\rm def}_L(\pi_K\eta^*)=\mbox{\rm def}_L(\eta^*)=b$
for $b=e_K+m$.
\end{lemma}
\begin{proof}
Recall that we have assumed that there is a $\eta\in L$ with norm
$N_{L/K}(\eta)=\eta^{\sigma+1}=(1+(\omega+\mu)^2\beta)(1+4\lambda)\equiv
1+(\omega^2+\mu^2)\beta\bmod \mu^2\beta\euP_K$. The congruence follows
from $m<b/2<e_K$.  We are interested in an explicit description for
$\eta\bmod \beta\euP_L$. So choose $\nu_0\in K$ with
$$v_K(\nu_0)=\begin{cases}
4e_K-2b&\mbox{for }b>4e_K/3,\\
2e_K-(b+1)/2&\mbox{for }b<4e_K/3.
\end{cases}$$
And observe that if $a\in 1+\nu_0\euP_K$ lies in the image of the norm
map $N_{L/K}$, we may assume that its preimage lies in $1+\beta\euP_L$
\cite[V\S3]{serre:local}. This means that we are really
interested in the image $\eta^{\sigma+1}\equiv 1+(\omega^2+\mu^2)\beta\bmod
(\mu^2\beta\pi_K,\nu_0\pi_K)$

Note that $(1+i\omega(x-1))^{\sigma+1}=1+\omega^2\beta-2i\omega$.
Choose a $2^f-1$ root of unity so that $(\omega')^2=\omega$.  Then
$(1+\omega'(i-1))^2=1-2i\omega+2\omega'(i-1)$ with
$v_K(2(i-1))=3e_K/2$.  By checking cases, we sees that
$v_K(2(i-1))>v_K(\nu_0)$ and that
$v_K(4)>v_K(2\beta)>v_K(\beta^2)>v_K(\nu_0)$. Therefore
$$[\eta\cdot (1+i\omega(x-1))\cdot (1+\omega'(i-1))]^{\sigma+1}\equiv
1+\mu^2\beta\bmod (\mu^2\beta\pi_K,\nu_0\pi_K).$$

Now observe that if $a'\in 1+\mu^2\beta\euO_K$ lies in the image of
the norm map $N_{L/K}$ with $\mbox{def}_K(a')=v_K(\mu^2\beta)$, we may
use \cite[V\S3]{serre:local} and choose its preimage $A$ to lie in
$1+B\euO_K$ with defect $\mbox{def}_L(A)=v_L(B)$ where
$$v_L(B)=\begin{cases}
2m+2e_K-b&\mbox{for }m+e_K<b,\\
4m+4e_K-3b&\mbox{for }m+e_K\geq b.
\end{cases}$$
As a result, $\eta\cdot (1+i\omega(x-1))\cdot (1+\omega'(i-1))\equiv
\eta^* \bmod (B\pi_L,\beta\pi_K)$ with
$\mbox{def}_L(\eta^*)=v_L(B)$. Since $v_L((i-1)(x-1))>v_L(\beta)$, we
can drop the ``$i$'' from $1+i\omega(x-1)$ and the first part of the
result follows.

Now consider the case $b=e_K+m$, which is equivalent to
$v_K(2/(\mu\beta))=0$. This means that there is a $q-1$ root of unity
$\omega_*$ such that $2/(\mu\beta)\equiv \omega_*\bmod
\euP_K$. Alternatively, this can be seen as the case when
$v_K(\mu^2\beta)=b$. So any element $a''\in1+\mu^2\beta\euP_K$ is a
norm from $L$ \cite[V\S3]{serre:local}.  This means that when
$\eta^*\not\equiv 1\bmod (B\pi_L,\beta\pi_K)$ (when it is relevant),
there is a $q-1$ root of unity $a$ such that
$\eta^*\equiv(1+a\mu(x-1))^{\sigma+1}\equiv 1+\mu^2\beta\bmod
\mu^2\beta\euP_K$, which means that the equation $a^2+a\omega_*\equiv
1\bmod\euP_K$ is solvable for $a$. Clearly
$a\not\in\{1,\omega_*^{-1}\}$.  Now note that because
$(x-1)^2=\beta\cdot(1-2(x-1)/\beta))$, we have
$\mu_K=\beta=(1-2(x-1)/\beta)$ in $L^*/(L^*)^2$. Thus $\pi_K\cdot
(1+\omega_*\mu(x-1))\in(L^*)^2(1+\euP_L^{b+1})$ and so
$\mbox{def}_L(\pi_K(1+\epsilon))=b=\mbox{def}_L(1+\epsilon)$.
\end{proof}

\begin{lemma}
If $i\in K$, $b> e_K$ and
$\omega^3\neq 1$, then
$$\mbox{\rm def}_{\mathcal{M}^*}(xy\eta)=\begin{cases}
\min\{4e_K-b+4m, 8e_K-3b\}&\mbox{\rm for }m+e_K<b,\\
\min\{8e_K-5b+8m, 8e_K-3b\} &\mbox{\rm for }m+e_K\geq b.
\end{cases}$$
\end{lemma}
\begin{proof}
Recall from (2) that $yY=1+(\omega+\mu)(x-1)$. Use Lemma 3.6 to find
that $x\cdot yY\cdot \eta = x\cdot (1+(\omega+\mu)(x-1))\cdot
(1+\omega(x-1))\cdot\eta^*$ in $\mathcal{N}^*$. Expand the product
$(1+(\omega+\mu)(x-1))\cdot (1+\omega(x-1))\equiv (1+\mu(x-1))\cdot
(1+\omega^2\beta)\bmod \beta\euP_L$.  So $x\cdot yY\cdot \eta = x\cdot
(1+\mu(x-1)) \cdot\eta^*$ in $\mathcal{N}^*$. Since
$\omega+\omega^2\equiv 1\bmod 2$, Lemma 3.4 yields
$Y=(1+(1+\mu+\mu^2)(x-1))\cdot (1+(x-1)(Y-1))$ in $\mathcal{M}^*$.
Thus $xy\eta=x\cdot (1+\mu(x-1)) \cdot\eta^*\cdot
(1+(1+\mu+\mu^2)(x-1))\cdot (1+(x-1)(Y-1))$ in $\mathcal{M}^*$.  In
general for $a\in\euO_K$, $(1+a(x-1))^2\equiv 1+a^2\beta\bmod
\beta\euP_L$ where of course $1+a^2\beta\in K$. Therefore
$xy\eta=\eta^*\cdot (1+\mu^2)(x-1))\cdot (1+(x-1)(Y-1))$ in
$\mathcal{M}^*$. Since $\mbox{def}_L(\eta^*)<v_L(\mu^2(x-1))$, we have
$\mbox{def}_L(\eta^*(1+\mu^2(x-1)))=\mbox{def}_L(\eta^*)$.  We need to
use Lemma 2.1 to determine $\mbox{def}_M(\eta^*)$. First note that
$\mbox{def}_L(\eta^*)<4e_K-b$. For $m+e_K<b$ this follows from the
fact that $m< e_K$ (otherwise $m\geq e_K$ and $b>m+e_K\geq 2e_K$,
yielding a contradiction). For $m+e_K\geq b$ this follows from
$m<b/2$.  Then because $\mbox{def}_L(\eta^*)<4e_K-b$, we have
$\mbox{def}_M(\eta^*)=2\mbox{def}_M(\eta^*)+b$. It is easy to check,
using the fact that $b$ is odd, that $\mbox{def}_M(\eta^*)\neq
v_L((x-1)(Y-1))$. And so we have determined that there is a $k_0$ such
that $\mbox{def}_M(k_0xy\eta)=\min\left \{\mbox{def}_M(\eta^*),
v_M((x-1)(Y-1))\right\}$. Thus, unless $b=e_K+m$, we can use Lemma 3.1
to find that $\mbox{def}_M(k_0xy\eta)\neq\mbox{def}_M(k)$ for all
$k\in K$ and the result follows.

When $b=e_K+m$, $\mbox{def}_M(\pi_K)=\mbox{def}_M(\eta^*)$ and we need
to be careful that $\mbox{def}_M(\pi_K\eta^*)$ does not exceed
$\mbox{def}_M(\eta^*)$. But this follows from the last part of Lemma 3.6.
\end{proof}

\subsubsection{\textbf{Most of unstable ramification: assume $b\leq e_K$.}}

Combining Lemma 3.1 with Lemma 3.9, we find that $3b\leq
\mbox{\rm def}_M(kxy\eta)\leq \min\{8e_K-5b+8m,8e_K-3b\}$.
Moreover, each of these possible values for $\mbox{def}_M(\alpha_k)$
is realized.

Again we start with a lemma that describes $\eta$.

\begin{lemma}
If $i=\sqrt{-1}\in K$, $b<e_K$ and $\omega^3=1$ there is in $\mathcal{L}^*$
the following coset identity
$$\eta= \left (1+\left
[\omega^2+\mu^2+\frac{\omega(\beta/2)^2+\omega^2(\beta/2)^3}{1+(\beta/2)^3}
\right ]\cdot(\beta/2)(x-1)\right) \cdot E$$
for some
$E\in 1+\mu^2(\beta/2)(x-1)\euP_L$.
\end{lemma}

\begin{proof}
Note that since $b<e_K$ we have $\mathbb{X}:=\beta/2\in\euP_K$.  We
are interested in an expression for $\eta\bmod
(\pi^2\mathbb{X}(x-1)\pi_L,\beta\pi_L)$. Now note that because
$b<e_K$, any element in $1+\mu^2\beta\euP_K$ has a preimage under the norm
$N_{L/K}$ that lies in $1+\mu^2\mathbb{X}(x-1)\euP_L$, and any element
of $1+\euP_K^{2e_K-(b-1)/2}$ has a preimage in $1+\beta\euP_L$
\cite[V\S3]{serre:local}. So we study the image
$N_{L/K}\eta=\eta^{\sigma+1}\bmod (\mu^2\beta\pi_K,\pi_K^{2e_K-\frac{b-1}{2}})$.

To compute this image, observe that for $M\in\euO_K$, $(1-
M\mathbb{X}(x-1))^{\sigma+1}=1+ (M-M^2\mathbb{X}^2)\beta$.  Now
set
$$M=\omega^2+\mu^2+\frac{\omega\mathbb{X}^2+\omega^2\mathbb{X}^3}{1+\mathbb{X}^3}.$$
Since $v_K(2)=e_K>b>2m=v_K(\mu^2)$ we can expand $M^2\bmod 2$ and find that
$$(1-
M\mathbb{X}(x-1))^{\sigma+1}\equiv 1+\left [\omega^2+\mu^2+\frac{\omega^2\mathbb{X}^3+\omega\mathbb{X}^5}{1+\mathbb{X}^6}\right ]\beta\bmod \mu^2\beta\euP_K.$$
Since $v_K(\beta^2)>2e_K-(b-1)/2$, we have
$(1-
M\mathbb{X}(x-1))^{\sigma+1}\equiv (1+(\omega^2+\mu^2)\beta)\cdot T$ with
$$T=
1+\frac{\omega^2\mathbb{X}^3+\omega\mathbb{X}^5}{1+\mathbb{X}^6}\beta\bmod
(\mu^2\beta\pi_K,\pi_K^{2e_K-\frac{b-1}{2}}).$$ Choose $\omega'$ to be
a $q-1$ root of unity such that $(\omega')^2=\omega$ and observe that
since $2\mathbb{X}=\beta$,
$$\left(1+\omega'(i-1)\frac{\omega^2\mathbb{X}^2+\mathbb{X}^3}{1+\mathbb{X}^3}\right
)^2\equiv T \bmod 2(i-1)\mathbb{X}^2\euO_K.$$ Since $e_K>b$ we have
$v_K(2(i-1)\mathbb{X}^2)\geq 2e_K-(b-1)/2$. So when $T\not\equiv
1\bmod (\mu^2\beta\pi_K,\pi_K^{2e_K-\frac{b-1}{2}})$, it lies in
$(1+\euP_K)^2$. Thus $\eta\in(1- M\mathbb{X}(x-1))(1+\euP_K)^2(1+\mu^2\mathbb{X}(x-1)\euP_L)(1+\beta\euP_L)$ and the result is proven.
\end{proof}

\begin{lemma}
If $i\in K$, $b<e_K$ and $\omega^3=1$ then $$\mbox{\rm
def}_{\mathcal{M}^*}(xy\eta)=\min\{4e_K-b+4m, 8e_K-3b\}.$$
\end{lemma}
\begin{proof}
We follow the proof of Lemma 3.5. In fact the first three paragraphs
(up through (5)) of that proof hold here verbatim.  So we begin at
the point where we substitute $X=\omega+\mu$ into (4), but examine
the result in $\mathcal{L}^*$ (instead of $\mathcal{L}$). Using the
fact that $v_L(2(x-1))>v_L(\beta)$ and
$v_L((1+i)(x-1)\beta/2)>v_L(\beta)$, we find that in
$\mathcal{L}^*$, we have the following coset identity:
$$1=1+(\omega^2+\mu^2)(x-1)\frac{1}{1+(\omega^2+\mu^2)\beta/2}+(1+i)(\omega+\omega^2+\mu+\mu^2)(x-1).$$
Note that for $a_1,a_2\in(x-1)\euO_L$ we have
$(1+a_1)(1+a_2)=1+a_1+a_2$ as cosets in $\mathcal{L}^*$. So recalling
(2), $yY=1+(\omega+\mu)(x-1)$, we find that there is an $E\in
1+\mu^2(\beta/2)(x-1)\euP_L$ so the the coset identity can be
rewritten as
$$yY=(1+i(1+\mu+\mu^2)(x-1))\cdot
\left (1+\frac{\omega(\beta/2)}{1+\omega^2(\beta/2)}(x-1)\right
)\cdot E.$$ Recall that $\omega^3=1$ and $\omega+\omega^2\equiv 1
\bmod 2$. Multiply both sides by $x$ and use (5) to see that
$xyY=(1+i(1+\mu+\mu^2)(x-1))\cdot (1+
\Theta_1(\beta/2)(x-1)
)\cdot E$ 
where
$$\Theta_1=\frac{1}{1+\beta/2}+\frac{\omega}{1+\omega^2(\beta/2)}.$$
Using Lemma 3.8, we multiply both sides by $\eta$.  This
results in the equivalence of cosets in $\mathcal{L}^*$,
$xyY\eta=(1+i(1+\mu+\mu^2)(x-1))\cdot \left
(1+\Theta_2(\beta/2)(x-1)\right )\cdot E'$, where $E'\in
1+\mu^2(\beta/2)(x-1)\euP_L$ and
$$\Theta_2=\omega^2+\mu^2+\frac{\omega(\beta/2)^2+
\omega^2(\beta/2)^3}{1+(\beta/2)^3}
+\frac{1}{1+\beta/2}+\frac{\omega}{1+\omega^2(\beta/2)}\equiv
\mu^2\bmod 2.$$ To see that $\Theta_2\equiv \mu^2\bmod 2$, expand the
geometric series and break the sum up into three partial sums
according to the equivalence class modulo $3$ of the exponents (on
the terms $(\beta/2)^n$). There are three partial sums: $n\equiv
0,1$, and $2\bmod 3$. Note that for $n\equiv 0,2\bmod 3$ the partial
sums are identically zero. For $n\equiv 1\bmod 3$ the partial sum is
$0\bmod 2$.

In summary, we have proven that there is an $E'\in
1+\mu^2(\beta/2)(x-1)\euP_L$ such that
$xyY\eta=(1+i(1+\mu+\mu^2)(x-1))\cdot \left
(1+\mu^2(\beta/2)(x-1)\right )\cdot E'\in\mathcal{L}^*$.
Now use Lemma 3.4 to see that we have the coset identity
$$xy\eta= \left(1+\mu^2(\beta/2)(x-1)\right)\cdot (1+(x-1)(Y-1))\cdot
E'\in\mathcal{M}^*.$$ Finally, by Lemma 2.1, since $m<b/2$ we have
$\mbox{def}_M((1+\mu^2(\beta/2)(x-1))\cdot E')=8m+8e_K-5b$. Hence
there is a $k_0\in K$ such that
$\mbox{def}_M(k_0xy\eta)=\min\{8m+8e_K-5b, 8e_K-3b\}$. Note that
since $b$ is odd, $8m+8e_K-5b\neq 8e_K-3b$. The result now follows
from Lemma 3.1.
\end{proof}

\section{Case 2: Two Break Biquadratic Extensions}

Let $M=K(\sqrt{u},\sqrt{v})$ as in \S2 and assume that $M$ embeds in a
quaternion extension. We are interested in determining
$\mbox{def}_M(\alpha_k)$.  Since the ramification filtration of
$\mbox{Gal}(M/K)$ is asymmetric with respect to the group action, we
have three cases to consider:
\begin{enumerate}
\item \textbf{$G_{b_2}$ fixes $K(\sqrt{u})$.} \ In this case, $b_1$ is the break
number of $K(\sqrt{u})/K$, $b_2$ is the break of $M/K(\sqrt{u})$,
$(b_1+b_2)/2$ is the break of $K(\sqrt{uv})/K$, and $b_1$ is the break
of $M/K(\sqrt{uv})$.
\item \textbf{$G_{b_2}$ fixes $K(\sqrt{uv})$.} \  In this case, $(b_1+b_2)/2$ is
the break of $K(\sqrt{u})/K$, $b_1$ is the break of $M/K(\sqrt{u})$
and $K(\sqrt{uv})/K$, and $b_2$ is the break of $M/K(\sqrt{uv})$.
\item \textbf{$G_{b_2}$ fixes $K(\sqrt{v})$.} \  In this case, $(b_1+b_2)/2$ is
the break of $K(\sqrt{u})/K$, $b_1$ is the break of $M/K(\sqrt{u})$,
$(b_1+b_2)/2$ is the break of $K(\sqrt{uv})/K$, and $b_1$ is the
break number of $M/K(\sqrt{uv})$.
\end{enumerate}

\subsection{Stable Ramification}

We begin by considering the case in which ramification is stable.  We
{\em do not} assume $\sqrt{-1}\in K$.  Because of the following lemma,
we conclude that if $b_1+b_2>2e_K$ then the third break number $b_3$
must be $4e_K+b_2$, which is precisely the value given in the catalog
in Section~\ref{ssec:catalog}.

\begin{lemma}
Suppose that  $b_1+b_2>2e_K$. Then
$\mbox{def}_M(\alpha_k)=4e_K-b_2$.
\end{lemma}
\begin{proof}
Because of the possibility that $\sqrt{-1}\not\in K$ we have
$\alpha_k=k\sqrt{uv}\eta\tau$.  The proof breaks naturally into
three steps. First we prove that $\mbox{def}_M(k)>4e_K-b_2$ for all
$k\in K$. We then show that $\mbox{def}_M(\tau)>4e_K-b_2$, and finally
prove that in each of the three cases
$\mbox{def}_M(\sqrt{uv}\eta)=4e_K-b_2$. The result will then follow
immediately.

We begin by considering $\mbox{def}_M(k)$. Choose $L$ to be the
fixed field of $G_{b_2}$, so that the break of $L/K$ is $b_1$. By
Lemma 2.1 we have $\mbox{def}_L(k)\geq b_1$. The break of $M/L$ is
$b_2$ and $b_1>2e_K-b_2$, and so by Lemma 2.1,
$\mbox{def}_M(k)>4e_K-b_2$.

For the second step, we need to consider $\mbox{def}_M(\tau)$.
Recall that $\tau\in K(\sqrt{uv})$, $N_{K(\sqrt{uv})/K}(\tau)=-1$
and note that $\mbox{def}_K(-1)\geq e_K$.  In Cases (1) and (3),
since $e_K<(b_1+b_2)/2$ using \cite[V\S3]{serre:local} we see that
$\mbox{def}_{K(\sqrt{uv})}(\tau)\geq e_K$.  Because
$e_K>2e_K-(b_1+b_2)/2>2e_K-b_1$, Lemma 2.1 with respect to
$M/K(\sqrt{uv})$ yields $\mbox{def}_M(\tau)>4e_K-b_1>4e_K-b_2$.  Now
consider Case (2): either $b_1>e_K$ or $b_1\leq e_K$. In the first
situation, using \cite[V\S3]{serre:local} we see that
$\mbox{def}_{K(\sqrt{uv})}(\tau)\geq e_K$. But since $b_2>b_1>e_K$,
$\mbox{def}_{K(\sqrt{uv})}(\tau)>2e_K-b_2$.  In the second
situation, using \cite[V\S3]{serre:local} we see that
$\mbox{def}_{K(\sqrt{uv})}(\tau)\geq 2e_K-b_1>2e_K-b_2$. Using Lemma
2.1 with respect to $M/K(\sqrt{uv})$ yields
$\mbox{def}_M(\tau)>4e_K-b_2$. This completes the second part of the
proof.

Finally,  we proceed to prove $\mbox{def}_M(\sqrt{uv}\eta)
=4e_K-b_2$. This time we will need to treat the three cases
separately.

Suppose first that we are in case (3). Since
$\mbox{def}_K(v)=2e_K-b_1$ to use \cite[V\S3]{serre:local}, we need
to consider the two possibilities $2e_K-b_1< (b_1+b_2)/2$ and
$2e_K-b_1\geq (b_1+b_2)/2$ separately. In both situations however we
see that $\mbox{def}_{K(\sqrt{u})}(\eta)>2e_K-(b_1+b_2)/2$,  so that
using Lemma 2.1 yields $\mbox{def}_M(\eta)>4e_K-b_2$. Note that
$\mbox{def}_{K(\sqrt{uv})}(\sqrt{uv})=2e_K-(b_1+b_2)/2$, which
implies that  $\mbox{def}_M(\sqrt{uv})=4e_K-b_2$. Hence we have
$\mbox{def}_M(\sqrt{uv}\eta)=4e_K-b_2$.

Case (2) is similarly easy. Note that in this case
$\mbox{def}_{K(\sqrt{uv})}(\sqrt{uv})=2e_K-b_1>2e_K-b_2$, and so
again by Lemma 2.1, $\mbox{def}_M(\sqrt{uv})>4e_K-b_2$. Since
$\mbox{def}_K(v)=2e_K-(b_1+b_2)/2$ and
$2e_K-(b_1+b_2)/2<(b_1+b_2)/2$, which is the break of
$K(\sqrt{u})/K$, we see that
$\mbox{def}_{K\sqrt{u})}(\eta)=2e_K-(b_1+b_2)/2$. Hence by Lemma 2.1
$\mbox{def}_M(\eta)=4e_K-b_2$ and so we have
$\mbox{def}_M(\sqrt{uv}\eta)=4e_K-b_2$.

In Case (1), $\mbox{def}_K(v)=2e_K-(b_1+b_2)/2$, and there are two
cases to consider depending upon whether $2e_K-(b_1+b_2)/2>b_1$ or
$2e_K-(b_1+b_2)/2\leq b_1$. In both cases, without loss of
generality we find that $\mbox{def}_{K(\sqrt{u})}(\eta)>2e_K-b_2$.
Therefore $\mbox{def}_M(\eta)>4e_K-b_2$. On the other hand,
$\mbox{def}_{K(\sqrt{uv})}(\sqrt{uv})=2e_K-(b_1+b_2)/2$, which means
that $\mbox{def}_M(\sqrt{uv})=4e_K-b_2$. Combining these yields
$\mbox{def}_M(\sqrt{uv}\eta)=4e_K-b_2$.
\end{proof}

\subsection{Unstable Ramification}
We assume here that $\sqrt{-1}\in K$.  Because of the following lemma,
if $b_1+b_2 < 2e_K$ then the third break number $b_3$ must be
$2b_1+3b_2$, $8e_K-2b_1-b_2$, or else must satisfy $2b_1+3b_2< b_3 <
8e_K-2b_1-b_2$ with $b_3 \equiv b_2$ mod $8$. These are precisely the
values listed in the catalog in Section~\ref{ssec:catalog}.
\begin{lemma}
Suppose that $M/K$ has two breaks which satisfy $b_1+b_2 < 2e_K$.
Then the value of $\mbox{def}_M(\alpha_k)$ must be $b_2+2b_1$,
$8e_K-3b_2-2b_1$ or else satisfy
$b_2+2b_1<\mbox{def}_M(\alpha_k)<8e_K-3b_2-2b_1$ with
$\mbox{def}_M(\alpha_k)\equiv -b_2\bmod 8$.
\end{lemma}
\begin{proof}
Because $\sqrt{-1}\in K$, $\alpha_k=k\sqrt{uv}\eta$.  Without loss of
generality we assume that $K(\sqrt{u})$ is the fixed field of
$G_{b_2}$. To proceed with the proof, we separate $\alpha_k$ into two
parts, $k\sqrt{uv}$ and $\eta$. Our first step is to explicitly
determine the possibilities for $\mbox{def}_M(k\sqrt{uv})$ from
$\mbox{def}_{K(\sqrt{uv})}(k\sqrt{uv})$. This will use the
classification of all possible second ramification numbers in a cyclic
degree $4$ extension from \cite{wyman}. Once that has been completed,
we will show that $\mbox{def}_M(k\sqrt{uv})<\mbox{def}_M(\eta)$. Since
this implies that $\mbox{def}_M(\alpha_k) = \mbox{def}_M(k\sqrt{uv})$,
the result will follow.

\textbf{Step 1.} Since $\sqrt{-1}\in K$, $K(\sqrt[4]{uv})/K$ is
cyclic of degree $4$.  Given $k\in K$ either
$K(\sqrt{k},\sqrt[4]{uv})/K$ is cyclic of degree $4$ or
$\mbox{Gal}(K(\sqrt{k},\sqrt[4]{uv})/K)\cong C_2\times C_4$. In
either case, the subextension $K(\sqrt{k}\sqrt[4]{uv})/K$ will be
cyclic of degree $4$. Again using the assumption that $\sqrt{-1}\in
K$, it is easy to check that any cyclic extension of degree $4$
containing $K(\sqrt{uv})$ is expressible as
$K(\sqrt{k}\sqrt[4]{uv})/K$ for some $k\in K$.

Let $t=(b_2+b_1)/2$. Then $t$ is the break number of
$K(\sqrt{uv})/K$. Since $t<e_K$, we find using \cite[Thm 32]{wyman}
that the second break number of $K(\sqrt{k}\sqrt[4]{uv})/K$ must be
one of  $b_2=3t$, $4e_K-t$ or $b_2=t+4m$ with $3t<b+4m<4e_K-t$, and
moreover that all these possible values for $b_2$ actually occur.
Therefore $\mbox{def}_{K(\sqrt{uv})}(k\sqrt{uv})$ must be one of the
values $t$, $4e_K-3t$ or
$t<\mbox{def}_{K(\sqrt{uv})}(k\sqrt{uv})<4e_K-3t$ with
$\mbox{def}_{K(\sqrt{uv})}(k\sqrt{uv})\equiv -t\bmod 4$.  Since in
every one of these cases
$\mbox{def}_{K(\sqrt{uv})}(k\sqrt{uv})<4e_K-b_1$, we may use Lemma
2.1 to determine that the possible values of
$\mbox{def}_M(k\sqrt{uv})$ are $b_2+2b_1$, $8e_K-3b_2-2b_1$ or
$b_2+2b_1<\mbox{def}_M(k\sqrt{uv})<8e_K-3b_2-2b_1$ with
$\mbox{def}_M(k\sqrt{uv})\equiv -b_2\bmod 8$.

\textbf{Step 2.} Recall that we have chosen $\eta\in K(\sqrt{u})$ so
that $N_{K(\sqrt{u})/K}(\eta)=v$. Since the break of $K(\sqrt{v})$
is $t=(b_2+b_1)/2$, the defect $\mbox{def}_K(v)=2e_K-t$. Moreover,
since $t<e_K$, we have $b_1<e_K$, and so $b_1<2e_K-t$. Since $b_1$
is the break number of $K(\sqrt{u})/K$, we use
\cite[V\S3]{serre:local} to choose $\eta\in K(\sqrt{u})$ such that
$\mbox{def}_{K(\sqrt{u})}(\eta)=4e_K-b_2-2b_1$, and since the break
of $M/K(\sqrt{u})$ is $b_2$ we may use Lemma 2.1 to find that
$\mbox{def}_M(\eta)=8e_K-b_2-4b_1$. To finish, it suffices to note
that $8e_K-b_2-4b_1> 8e_K-3b_2-2b_1\geq \mbox{def}_M(k\sqrt{uv})$.
\end{proof}

\section{Proof of Main Results}
\begin{proof}[Proof of Theorem 1.1]
Given the computations in Sections 3 and 4, the first statement has
already been established; we just need to prove the second. In fact,
since $\sqrt{-1}\in K$ it is enough to prove, for each
$i\in\{1,1^*,2\}$, that if $s_1\leq s_2$ are the first two coordinates
of a triple in $\mathcal{R}_i^e$ then there exist a pair of elements
$u,v\in K^*$ such that $(u,v)=1$ and the following:
\begin{enumerate}
 \item
If $i=2$, then
$v_K(u-1)=2e_K-s_1$ and $\mbox{def}_K(v)=2e_K-(s_2+s_1)/2$.  
\item
If $i\in\{1,1^*\}$, then $v_K(u-1)=v_K(v-1)=2e_K-s_1$. Indeed 
$v-1=
(\omega+\mu)^2(u-1)$ where
\begin{enumerate}
\item
$\omega$ is a nontrivial $q-1$ root of unity, with
$\omega^2+\omega+1=0$ if $i=1^*$ and $\omega^3\neq 1$ if $i=1$, and
\item $\mu\in K$, with 
$v_K(\mu)=(s_2-s_1)/4$ if $s_2<\min\{2s_1,4e_K-s_1\}$ and
$v_K(\mu)>(s_2-s_1)/4$ if $s_2=\min\{2s_1,4e_K-s_1\}$.
\end{enumerate}
\end{enumerate}
Under these circumstances, the biquadratic extension
$M=K(\sqrt{u},\sqrt{v})$ will embed in a quaternion extension, which
lies in $\mathcal{Q}_i^K$ and have the desired ramification pair
$s_1\leq s_2$. The third coordinate $s_3$ of the triple is then
achieved by using Lemma 3.1 to choose the $k$ in $\alpha_k$
appropriately.

We begin with cases where any pair of elements $u,v$ with the desired
defects must automatically satisfy $(u,v)=1$.
If
$i\in\{2\}$ and $s_2+3s_1<4e_K$, then every pair of elements $u,v\in
K^*$ that satisfies $v_K(u-1)=2e_K-s_1$ and
$\mbox{def}_K(v)=2e_K-(s_2+s_1)/2$ will satisfy $(u,v)=1$
\cite[V\S3]{serre:local}.  Similarly for $i\in\{1,1^*\}$ when
$s_1<e_K$ every pair of units $u,v\in 1+\euP_K$ that satisfies
$v_K(u-1)=v_K(v-1)=2e_K-s_1$ will satisfy $(u,v)=1$. In these cases
every such biquadratic extension $K(\sqrt{u},\sqrt{v})/K$ embeds.

Outside of these two cases, we are free to choose $u$ based upon
defect alone, but must choose $v$ dependent upon $u$.  Suppose
$i\in\{2\}$ and $s_2+3s_1\geq 4e_K$. Pick any element $u\in K^*$ such
that $v_K(u-1)=2e_K-s_1$. Pick any element $\nu\in K(\sqrt{u})$ with
$v_{K(\sqrt{u})}(\nu-1)=2e_K-(s_2+s_1)/2$. By
\cite[V\S3]{serre:local} $N_{K(\sqrt{u})/K}(\nu)=v\in K^*$ where
$v_K(v-1)=2e_K-(s_2+s_1)/2$ and by design, $(u,v)=1$.

Consider now $i\in\{1,1^*\}$ and $s_1> e_K$.  Pick $\omega$ according
to whether $i=1$ or $i=1^*$. Choose $\beta\in K$ with
$v_K(\beta)=2e_K-s_1$. Let $x^2=u=1+\beta$. Then $L/K$ where $L=K(x)$
has ramification break number $s_1$. We consider the cases
$s_1<3e_K/2$ and $s_1\geq 3e_K/2$ separately. 

If $e_K<s_1<3e_K/2$,
then then one can show, as in the proof of Lemma 3.6, that the norm
$$N_{L/K}\bigl((1+i\omega(x-1))\cdot (1+\omega'(i-1))\bigr)\equiv
1+\omega^2\beta\bmod \euP_K^{s_1+1}$$ where $(\omega')^2=\omega$.  By
\cite[V\S3]{serre:local}, any element $\alpha\in 1+\euP_K^{s_1+1}$
satisfies $(1+\beta,\alpha)=1$. As a result,
$(1+\beta,1+\omega^2\beta)=1$. Again using \cite[V\S3]{serre:local},
there is, for any relevant value of $m$, an element $A\in L$ such that
$v_K(N_{L/K}(A)-1)=2e_K-b+2m$. Note $(1+\beta,(1+\omega^2\beta)\cdot
N_{L/K}(A))=1$. Using 
Lemma 2.2, we find a $\mu\in K$
with $v_K(\mu)=m$ such that $1+(\omega+\mu)^2\beta\equiv
(1+\omega^2\beta)\cdot N_{L/K}(A)$ (modulo squares in $K^*$).
So for $v=1+(\omega+\mu)^2\beta$ we have $(u,v)=1$.

Assume $s_1\geq 3e_K/2$.
Then following \cite[V\S3]{serre:local} we see explicitly
that for $a\in \euO_T$,
$$N_{L/K}\left (1+a\frac{2}{1-x}\right)=1-\frac{4}{\beta}(a+a^2).$$
Let $\lambda\in \euO_T$ such that $y+y^2=\lambda$ is irreducible. Then
$(1+\beta,1-4\lambda/\beta)=-1$. This means that we have
$(1+\beta,1+\omega^2\beta)=1$ or
$(1+\beta,1+\omega^2\beta-4\lambda/\beta)=1$.  If we express
$1+\omega^2\beta-4\lambda/\beta$ as $(1+(\omega+\mu_0)^2\beta)\bmod 4$
using Lemma 2.2 we find that $v_K(\mu_0)=m=s_1-e_K$. Of course,
$1+\omega^2\beta=(1+(\omega+\mu_0)^2\beta)$ with $v_K(\mu_0)=\infty>
s_1-e_K$. Let $v_0= 1+(\omega+\mu_0)^2\beta$. Then $(u,v_0)=1$ and
$v_K(\mu_0)=m\geq s_1-e_K$.  Since $s_1+4m>4e_K-s_1$ we find that
$K(\sqrt{u},\sqrt{v_0})$ embeds and has refined ramification
filtration $s_1<s_2=4e_K-s_1=\min\{2s_1,4e_K-s_1\}$.  All that remains
is the situation where $s_2<4e_K-s_1=\min\{2s_1,4e_K-s_1\}$.  In other
words, for each $0<m<e_K-b/2$, we must find $\mu\in K$ with
$v_K(\mu)=m$ so that $v=1+(\omega+\mu)^2\beta$ satisfies
$(u,v)=1$. The second refined break of $K(\sqrt{u},\sqrt{v})$ will
then be $s_2=s_1+4m<4e_K-s_1=\min\{2s_1,4e_K-s_1\}$.  We know that
$s_1\geq 3e_K/2>4e/3$. We have $2m<2e_K-s_1$.  Pick any $A\in L$ with
$v_L(A)=2m+2e_K-s_1<4e_K-2s_1<s_1$.  Then using
\cite[V\S3]{serre:local}, $v_K(N_{L/K}(A)-1)=2e_K-b+2m$. So
$(u,v_0\cdot N_{L/K}(A))=1$. And using Lemma 2.2, we can express
$v_0\cdot N_{L/K}(A)$ as $v\equiv 1+(\omega+\mu)^2\beta\bmod 4$
(also modulo squares in $K^*$) where $\mu\in K$ with $v_K(\mu)=m$.
\end{proof}
\begin{proof}[Proof of Theorem 1.2]
This follows immediately from \S3.3.1 and \S4.1.
\end{proof}

\end{document}